\documentclass[twocolumn,superscriptaddress,longbibliography,floatfix,aps,pre,showkeys]{revtex4-2}

\usepackage{graphicx}
\usepackage[dvipsnames]{xcolor}
\usepackage{amsbsy}
\usepackage{tikz}
\usepackage{amsthm,mathrsfs,amsopn,amsfonts,amsmath,amssymb}
\usepackage[justification=RaggedRight]{caption}

\usepackage{comment}
\usepackage{etoolbox}
\AtBeginEnvironment{matrix}{\everymath{\displaystyle}}
\usepackage{url}
\usepackage{multirow}

\newcommand{\rev}{\color{black}}

\usepackage[colorlinks=true,linkcolor=blue,citecolor=blue,urlcolor=blue]{hyperref}

\usepackage[capitalise]{cleveref}

\providecommand{\keywords}[1]
{
  \small	
  \textbf{Keywords--} #1
}

\begin{document}

\title{Hamiltonian control to desynchronize Kuramoto oscillators with higher-order interactions}

\author{Martin Moriamé}
\email{martin.moriame@unamur.be} 
\affiliation{Université de Namur (Belgium), Department of Mathematics \& Namur Institute for Complex Systems, naXys}

\author{Maxime Lucas}
\affiliation{Université de Namur (Belgium), Department of Mathematics \& Namur Institute for Complex Systems, naXys}
\affiliation{Université Catholique de Louvain, Earth and Life Institute, Mycology, B-1348, Louvain-la-Neuve, Belgium}

\author{Timoteo Carletti}
\affiliation{Université de Namur (Belgium), Department of Mathematics \& Namur Institute for Complex Systems, naXys}

\date{\today}

\begin{abstract}
	Synchronization is a ubiquitous phenomenon in nature. Although it is necessary for the functioning of many systems, too much synchronization can also be detrimental, e.g., (partially) synchronized brain patterns support high-level cognitive processes and bodily control, but hypersynchronization can lead to epileptic seizures and tremors, as in neurodegenerative conditions such as Parkinson's disease.
    Consequently, a critical research question is how to develop effective pinning control methods capable to reduce or modulate synchronization as needed.
    Although such methods exist to control pairwise-coupled oscillators, there are none for higher-order interactions, despite the increasing evidence of their relevant role in brain dynamics. 
    In this work, we fill this gap by proposing a generalized control method designed to desynchronize Kuramoto oscillators connected through higher-order interactions. Our method embeds a higher-order Kuramoto model into a suitable Hamiltonian flow, and builds up on previous work in Hamiltonian control theory to analytically construct a feedback control mechanism.
    We numerically show that the proposed method effectively prevents synchronization {\rev in synthetic and empirical higher-order networks}. Although our findings indicate that pairwise contributions in the feedback loop are often sufficient, the higher-order generalization becomes crucial when pairwise coupling is weak. Finally, we explore the minimum number of controlled nodes required to fully desynchronize oscillators coupled via an all-to-all hypergraphs.
\end{abstract}


\maketitle


\section{Introduction}

Synchronization, the emergence of order in the collective dynamics of coupled oscillators, is key to many natural and man-made systems~\cite{pikovskysynchronization,strogatz2004}. Synchronization can be found in domains ranging from physics to biology and neuroscience, with typical examples including the clapping in unison of an audience after a concert, or the (hyper)synchronized firing of neurons in the brain at the onset of an epileptic seizure~\cite{cumin2007generalising}. The paradigmatic model of synchronization is that of phase oscillators  all-to-all coupled, first introduced by Y. Kuramoto forty years ago~\cite{kuramoto1984chemical}, which led to breakthroughs in our understanding of collective dynamics. Since then, {\rev the} Kuramoto model has been extended in many ways, most notably to incorporate a complex network of (pairwise) interactions between oscillators~\cite{arenas2008,rodrigues2016}. 

More recently, increasing experimental and theoretical evidence suggests that  networks, or pairwise coupling schemes, may be not precise enough when modeling synchronization and complex systems in a more general setting. In fact, these systems often appear to be better modeled by higher-order (i.e., group) interactions between any number of oscillators at a given time~\cite{battiston2020networks,battiston2021physics}. These higher-order interactions have been shown to dramatically affect dynamics in many processes such as consensus~\cite{neuhauser2020multibody,deville2021consensus}, spreading~\cite{iacopini2019simplicial,lucas2023simplicially,ferrazdearruda2023multistability}, diffusion~\cite{carletti2020random,schaub2018flow}, and evolution~\cite{alvarez-rodriguez2021evolutionary}. In coupled oscillators, these higher-order interactions have been shown to naturally appear from phase reduction of pairwise coupled nonlinear oscillators~\cite{bick2016chaos,gengel2020highorder,nijholt2022emergent}. They have also been shown to favor explosive transitions~\cite{skardal2019abrupt,kuehn2020universal}, chaos~\cite{bick2016chaos}, multistability~\cite{tanaka2011multistable, matheny2019exotic,zhang2023deeper}, and chimera-type states~\cite{kundu2022higherorder,MNGCF2024} and more~\cite{nurisso2023unified}.

In many cases, synchronization can be beneficial, for example, in adjusting human circadian rhythms to the day-night cycles~\cite{Hafner2012}. Synchronization can, however, also be detrimental, as in epileptic seizures~\cite{cumin2007generalising}. The latter being characterized by a hyper-synchronization state that is abnormally strong, persistent in time and acting on a large portion of the brain, whose ultimate effect is to induce malfunctioning in the patient behavior. Methods to precisely control the system to avoid or reduce synchronization can thus have a crucial impact on the systems, and in the case of epilepsy, on the well-being of the patients. {\rev In this context, pinning control methods, i.e. control signals injected in particular locations of the system, are well suited~\cite{tass,neurocontrol,DBS}. 

Pinning control methods are of great interest in the context of complex systems; they have been widely studied in complex networks~\cite{wang2002pinning,Wenwu2009onpinning,liu2011controllability,liu2016control,liu2021optimizing,Moriame2023onthelocation} and have recently been extended to higher-order networks~\cite{Chen2021controllability,delellis2023pinninglinear,dellarossa2023emergence,rizzello2024pinning,Muolo2024Pinning,wang2024synchronization,Li2024synchronization,Xia2024pinning}. For instance, authors in Refs.~\cite{delellis2023pinning,delellis2023pinninglinear,dellarossa2023emergence,rizzello2024pinning} found that the linearized model on directed higher-order networks could be formulated in an elegant way by using signed a network, i.e. a network whose edges can have positive or negative weights. They then established synchronizability criteria from this formulation.
Let us stress that these studies focused on the stabilization of the synchronized state; to the best of our knowledge, methods designed to reduce  synchronization do exist solely in the case of pairwise-coupled oscillators~\cite{gjata2017using,asllani2018minimally}. The aim of this work is to fill this gap by proposing a control method allowing to reduce synchronization for higher-order networks. This result is thus especially important, because higher-order interactions can favor explosive transitions to synchronization~\cite{skardal2019abrupt, millan2020explosive, landry2020effect,st-onge2022influential, burgio2024triadic, malizia2025hyperedge} and induce stronger local stability of the synchronized state~\cite{zhang2023deeper}.}

In this paper, we thus propose a control method to desynchronize oscillators coupled via higher-order interactions, by generalizing the method developed in~\cite{gjata2017using,asllani2018minimally} for pairwise interactions. {\rev To derive this control, we (i) found a Hamiltonian system admitting an invariant torus whose dynamics return the one of the higher-order Kuramoto model (HOKM), and (ii) we generalized Hamiltonian control techniques ~\cite{vittot2004perturbation,ciraolo2004control} to control and desynchronize oscillators coupled through higher-order interactions. Finally, we demonstrate the efficiency of the proposed control method on both synthetic and empirical higher-order networks. Moreover the control term is adaptive, namely it is always present, i.e., we are not considering switched-systems, it reaches large values once the original uncontrolled system naturally synchronizes, but it remains small when the system spontaneously evolves in a disordered state.}


The derived control method contains thus pairwise terms as well as many-body ones (we restrict our analysis to $3$-body terms in this work, but this assumption can be clearly relaxed). {\rev In the presence of pairwise interactions alone, we recover the results of~\cite{gjata2017using}. We then numerically show that the control method is effective when there are either pairwise terms, $3$-body interactions or both of them, i.e., it manages to significantly reduce synchronization of the system. Furthermore, we compare it to the pairwise version~\cite{gjata2017using} and show that in a quite large number of cases the latter is sufficient, with the notable exception of the case where the pairwise coupling is small compared to the higher-order one, in this case the higher-order generalization is required to achieve desynchronization.} Finally, we show that synchronization can be impeded even when the control acts only on a fraction of nodes, we can thus conclude that the proposed control method is not very invasive. Our results suggest the need for a critical proportion of controlled nodes to drastically decrease synchronization in all-to-all hypergraphs.

The paper is organized as follows. In Sec. \ref{sec:HOKM}, we define a general HOKM, then in Sec. \ref{sec:methods}, we present the general formulation of the Hamiltonian system that embeds the HOKM {\rev and we analytically derive the method to control the HOKM .} The numerical validation is presented in Sec. \ref{sec:results}. Finally, we discuss the results and conclude in Sec \ref{conclusion}.

\section{The Higher-Order Kuramoto model}
\label{sec:HOKM}

In this work, we focus on the following Higher-Order version of the Kuramoto model (HOKM) composed by $N$ non-identical oscillators, described by the angular variables $\theta_i$, $i=1,\dots,N$, interacting through first and second order (i.e. $2$- and $3$-body) interactions
\begin{equation}\label{eq:HOKM3}
    \begin{split}
        \dot{\theta_i}&=\omega_i+\frac{K_1}{N}\sum_{j=1}^NA_{ij}\sin(\theta_j-\theta_i)+\\
        &+\frac{K_2}{N^2}\sum_{j,k=1}^N B_{ijk}\left[\sin(\theta_j+\theta_k-2\theta_i)+\sin(2\theta_j-\theta_k-\theta_i)\right]\,
    \end{split}
\end{equation}
where $\mathbf{A}$ and $\mathbf{B}$ are the first- and second-order adjacency tensors that encode the interactions: $A_{ij} = 1$ if there is a link connecting oscillators $i$ and $j$, namely a first-order (pairwise) interaction between nodes $i$ and $j$, and $A_{ij} = 0$ otherwise. Similarly, $B_{ijk} = 1$ if the oscillators $i$, $j$ and $k$ interact together, namely a second-order (i.e., $3$-body) interaction between them, and $B_{ijk} = 0$ otherwise. For the second-order interactions,
we require $i \neq j \neq k$, so that each triplet involves three distinct nodes. The parameters $K_1\geq 0$ and $K_2\geq 0$ are, respectively, the first- and second-order coupling strengths, and $\omega_i$ is the natural frequency of oscillator $i$. 

Note that we consider undirected hypergraphs so that $\mathbf{A}$ is symmetric, $A_{ij}=A_{ji}$, and  $\mathbf{B}$ is also invariant under indices permutations, i.e., $B_{\pi(ijk)}$ keeps the same value for all permutations $\pi(ijk)$ of the three indices. 
In addition, we assume the natural frequencies $\omega_i$ to be non-resonant, namely $\forall \mathbf{k}\in\mathbb{Z}^N: \mathbf{k} \cdot (\omega_1,\dots,\omega_N)^\top =0$ if and only if $\mathbf{k}=0$. This assumption{\rev, although necessary in the following,} is not too restrictive since, as often assumed in the literature, $\omega_i$ are sampled from a continuous symmetric distribution and one can prove that the set of resonant frequencies has in this case zero measure.
The parameters $K_1$ and $K_2$ are normalized, respectively, by $N$ and $N^2$, the respective numbers of terms in each sum, to fairly compare structures of different sizes~\cite{arenas2008} and different amounts of first- and second-order interactions.

When $K_2=0$, \cref{eq:HOKM3} recovers the canonical Kuramoto model~\cite{kuramoto1984chemical} for complex networks. 
For $K_2>0$, the second term of \cref{eq:HOKM3} encodes the second-order interactions between oscillators as a combination of the two distinct second-order coupling functions 
\begin{equation}\label{3orderTerm1}
         C_1=\sin(\theta_j+\theta_k-2\theta_i)\, ,
\end{equation}
and
\begin{equation}\label{3orderTerm2}
         C_2=\sin(2\theta_k-\theta_j-\theta_i)\, .
\end{equation}
%

Although most of the literature considers either of these coupling functions~\cite{skardal2019abrupt,lucas2020multiorder,tanaka2011multistable}, systems derived from phase reduction approaches often display a weighted combination of both~\cite{bick2016chaos, leon2019phase,gengel2020highorder}. We know that they induce a difference in the speed of convergence to full synchronization~\cite{lucas2020multiorder} and may have other effects on the dynamics, but to date, there is no consensus on the best way to model second-order interactions in general. 
However, close to the synchronization manifold these two coupling functions only differ in the speed of convergence toward synchronization: it is twice faster with the coupling in \cref{3orderTerm1}, as one can prove with linear stability analysis~\cite{lucas2020multiorder}. It is an open question to study their impact on  desynchronization.
Interestingly, the combination of coupling functions from \cref{3orderTerm1} and \cref{3orderTerm2} arising in the HOKM under consideration \cref{eq:HOKM3} naturally derives from the Hamiltonian embedding system we present in the next section.


Finally, let us observe that the proposed model \cref{eq:HOKM3} can be defined more generally for interactions of any order $d> 2$, as we will see in Sec. \ref{sec:Ham}. However, the number of possible interaction terms grows exponentially with $d$ and thus the computations become more cumbersome. For the sake of clarity, we restrict our analysis to the cases $d=1$ and $d=2$.

{\rev 
\section{Hamiltonian Control Theory}
\label{sec:methods}
}
{\rev The aim of this section is to introduce the basis of Hamiltonian control theory adapted to model under investigation. We refer the interested reader to~\cite{vittot2004perturbation,ciraolo2004control} for a more general and detailed presentation.}

{\rev
\subsection{Hamiltonian system embedding HOKM}
\label{sec:Ham}
}

In~\cite{witthaut2014kuramoto}, authors proposed the Hamiltonian function
\begin{equation}\label{def:Ham}
	H=\sum_{i=1}^N \omega_i I_i - \frac{K_1}{N}\sum_{i,j=1}^N A_{ij}\sqrt{I_iI_j}(I_j-I_i)\sin(\theta_j-\theta_i)
\end{equation}
{\rev where the action variables, $I_i$, are the momentum variables, canonically conjugated to the angles, $\theta_i$}. This Hamiltonian is such that, for all positive constant $c$, the torus $T_c:=\{(\mathbf{I},\pmb{\theta})\in\mathbb{R}_+^N\times[0,2\pi]^N | \forall i\in\{1,\dots,N\}:I_i=c\}$ is invariant by the flow. Moreover, the solutions of the Hamiltonian system restricted to the torus $T_{\frac{1}{2}}$ exhibit angle variables evolving according to the classical Kuramoto model with natural frequencies $(\omega_1,\dots,\omega_N)$, coupling strength $K_1$ and coupling network adjacency matrix $\mathbf{A}$. The restriction to $T_c$ with $c\neq\frac{1}{2}$ also gives rise to the KM but with a coupling strength equal to $2cK_1$.

We hereby propose a straightforward generalization of the above Hamiltonian function given by
\begin{widetext}
    \begin{eqnarray}
    \label{def:HOHam3}
		H(\mathbf{I},\pmb{\theta})&:=&\sum_{i=1}^N I_i\omega_i - \frac{K_1}{N}\sum_{i,j=1}^N A_{ij}\sqrt{I_iI_j}(I_j-I_i)\sin(\theta_j-\theta_i) - \frac{K_2}{N^2}\sum_{i,j,k=1}^N B_{ijk}\sqrt[3]{I_iI_jI_k}(I_j+I_k-2I_i)\sin(\theta_j+\theta_k-2\theta_i)\notag\\
        &&\\
        &=&H_0(\mathbf{I})+V(\mathbf{I},\pmb{\theta})\, ,\notag
    \end{eqnarray}
\end{widetext}
where we added a new term involving three action variables and three angles variables, thus encoding for the three-body interactions. {\rev The Hamiltonian system defined by \eqref{def:HOHam3} preserves again the torus $T_c$ and the angle dynamics restricted on the latter torus, coincides with~\cref{eq:HOKM3} with $c=\frac{1}{2}$ (see \cref{app:Embeding}).}

Let us also observe {\rev (see \cref{app:Embeding})} that the chosen Hamiltonian system~\eqref{def:HOHam3} {\rev returns a control strategy involving both two three-body interaction terms~\eqref{3orderTerm1} and~\eqref{3orderTerm2}, as they derive from the triple sum term. This motivated our choice to work with the model \cref{eq:HOKM3} where the two possible $3$-body interactions are allowed. Let us observe that if we were interested in controlling an HOKM involving only one kind of $3$-body interaction, we could first compute the full control term and then remove the undesired terms so to get an effective control term.} 

Moreover, the proposed higher-order Hamiltonian system exhibits another important property, similarly to~\cref{def:Ham}:  angles synchronization induces an instability of $T_c$ in the transversal directions, i.e., orthogonal to $(\mathbf{I},\pmb{\theta})=(c\mathbf{1},\pmb{\theta})$ {\rev (see \cref{app:Embeding} for more details)}.

\vspace{0.5cm}

{\rev \subsection{Construction of the Hamiltonian control}
\label{sec:control}
}

The goal of this work is to design a feedback control, namely a term depending on the angles, $\pmb{\theta}$ that, once added to~\cref{eq:HOKM3} would prevent the system from synchronizing. Because of the link presented in the previous section between synchronizability and instability of the invariant torus, we are thus looking for a control term that should increase the stability of the invariant torus.

To achieve this, we make use of the Hamiltonian control theory developed in~\cite{vittot2004perturbation,ciraolo2004control}. Let us first remark that the Hamiltonian system $H(\mathbf{I},\pmb{\theta})$ given by Eq.~\eqref{def:HOHam3} can be decomposed as a sum of two terms, $H:=H_0(\mathbf{I})+V(\mathbf{I},\pmb{\theta})$, where $H_0(\mathbf{I}):=\sum_{i=1}^N\omega_iI_i${\rev .} 
As the coefficients $\frac{K_1}{N}$ and $\frac{K_2}{N^2}$ are in general smaller than $1$, the term $V:=H-H_0$ can be considered to be a perturbation of the integrable part $H_0$. 

 {\rev Let us briefly explain the main idea of Hamiltonian control theory and thus the rationale behind the proposed control strategy.}
 The theory of Hamiltonian control~\cite{vittot2004perturbation,ciraolo2004control} is based on the possibility to build a small but not null perturbation, $f(V)$, that once added to the Hamiltonian $H$ acts as a feedback control allowing the flow induced by $H_0+V+f(V)$ to be canonically conjugate to the one induced by $H_0$. This means, roughly speaking, that the controlled system behaves like a set of uncoupled oscillators described by $H_0$, each evolving independently from the others, with incommensurable frequencies, and thus no synchronization is possible.

{\rev The general expression of such term $f(V)$ is proposed in~\cite{vittot2004perturbation,ciraolo2004control} and can be formally written as}
\begin{equation} \label{eq:h_i}
		h_i^{(N)}:=-\frac{1}{2}\frac{\partial }{\partial I_i}\left[\{\Gamma V\}V\right]\big|_{\mathbf{I}=\frac{1}{2}}\, ,
\end{equation}
{\rev  where $\Gamma$ is the pseudo-inverse operator associate to the unperturbed Hamiltonian $H_0$ (the interested reader can find more details in~\cref{App:GeneralControlTerm})}. Note that as we evaluate the derivative on ${T}_{1/2}$, $h^{(N)}_i$ no longer depends on $\mathbf{I}$ but only on $\pmb{\theta}$ and the parameters of the system (i.e., $\mathbf{A}$, $\mathbf{B}$ and $\pmb{\omega}$). {\rev As a matter of notation, let us observe that} the exponent $N$ indicates that the control term is applied to all $N$ nodes.

Note that the control term, {\rev whose formal expression is given by \cref{eq:h_i}, is composed by several terms (see~\cref{App:GeneralControlTerm}) and thus it can be difficult to understand the role of the involved terms and not to be easily implementable in real applications. Moreover} it depends on all angles, while {\rev classically} (feedback) pinning control {\rev methods rely on the action on} a small number of nodes to control a complex (higher-order) network~\cite{wang2002pinning,liu2011controllability,liu2016control,asllani2018minimally,liu2021optimizing}. Let us therefore consider a case where a subset of $M\leq N$ nodes, that we can label without loss of generality as nodes $1,\dots,M$, receive {\rev the controlled signal}. In other words, one would like to be able to force the system to desynchronize by only acting on the dynamics of those nodes and letting the others evolve according to~\cref{eq:HOKM3}. Moreover, the injected control term should be computed only from the observed dynamics of those same nodes, as is classically done in feedback-loop pinning control.

Let us consider the following modification of $V$ defined in~\eqref{def:HOHam3} where we only take into account the $M$ first nodes, that is
\begin{eqnarray}\label{eq:split_V}
		V^{(M)}&=&
        V^{(M,1)}+V^{(M,2)}\\
        &=&-\frac{K_1}{N}\sum_{i,j=1}^MA_{ij}\sqrt{I_iI_j}(I_j-I_i)\sin(\theta_j-\theta_i) + \nonumber\\
		&-& \frac{K_2}{N^2}\sum_{i,j,k=1}^M B_{ijk}\sqrt[3]{I_iI_jI_k}(I_j+I_k-2I_i)\times \nonumber\\
        &&\times\sin(\theta_j+\theta_k-2\theta_i)\, , \nonumber
\end{eqnarray}
and then by following a similar analysis of the one above presented we can {\rev compute the control term by using}
\begin{equation}\label{eq:h_iPartial}
	h^{(M)}_i:=
		-\frac{1}{2}\frac{\partial }{\partial I_i}\left[\{\Gamma V^{(M)}\}V^{(M)}\right]\big|_{\mathbf{I}=\frac{1}{2}} {\rev \chi_i} \, .
\end{equation}
{\rev where $\chi_i=0$ if $i\in\{1,\dots,M\}$ and 0 otherwise.} Then  by adding~\cref{eq:h_iPartial} to~\cref{eq:HOKM3} {\rev we obtain} a novel pinning control scheme {\rev depending} an \textit{a priori} chosen subset of $M$ nodes.

Then, we can further remark that computing the control terms can be very costly (see Appendix \ref{App:FullControlTerm}) in particular due to the term arising from $V^{(M,2)}$, namely the $3$-body interaction. Let us then define a simplified control involving only $V^{(M,1)}$, i.e., by only using the contributions of  pairwise interactions and still restricted on the use of $M$ arbitrarily chosen nodes, {\rev once again the formal expression of the control can be obtained via the formula}
\begin{equation}\label{eq:h_iPartial2D}
	\tilde{h}^{(M)}_i:=
		-\frac{1}{2}\frac{\partial }{\partial I_i}\left[\{\Gamma V^{(M,1)}\}V^{(M,1)}\right]\big|_{\mathbf{I}=\frac{1}{2}} {\rev \chi_i} \, .
\end{equation}
Let us observe that the latter coincide with the {\rev control} studied in~\cite{asllani2018minimally}. {\rev The interested reader is referred to~\cref{App:FullControlTerm} for the explicit computations returning the control terms starting from  $h_i^{(M)}$ and $\tilde{h}_i^{(M)}$.}

In conclusion the equations ruling the dynamics of the controlled higher-order Kuramoto model are given by
\begin{equation}\label{eq:ControlledHOKM3}
	\begin{split}
		\dot{\theta_i}=&~\omega_i+\frac{K_1}{N}\sum_jA_{ij}\sin(\theta_j-\theta_i)+\\
		+&\frac{K_2}{N^2}\sum_{j,k} B_{ijk}\left[\sin(\theta_j+\theta_k-2\theta_i)-\sin(2\theta_j-\theta_k-\theta_i)\right] \\
        &+ h_i 
	\end{split}
\end{equation}
where $h_i=h_i^{(M)}$, {\rev in the case of the whole control term acting on $M$ nodes} or $\tilde{h}_i^{(M)}$ {\rev for the control limited to pairwise terms, still acting on $M$ nodes}.

{\rev 
A complete and rigorous interpretation of the impact of each term in the control function, is a difficult task. We can nevertheless gain some insights by adapting to the present framework the analysis proposed for the pairwise case in Ref.~\cite{asllani2018minimally}. There, authors shown that, in the case of complete graphs, the dominant term is given by 
\begin{equation}\label{eq:pairwise_reduction}
    \tilde{h}_i^{(N)}\approx -\frac12 K_1^2R\tilde{R}_i\cos(\Psi-\tilde{\Psi}_i)\, ,
\end{equation}
where
\begin{equation}\label{eq:order_parameter}
    Re^{\iota\Psi} = \frac{1}{N} \sum_j e^{\iota\theta_j} \, ,
\end{equation}
is the Kuramoto order parameter, being $\iota:=\sqrt{-1}$, and
\begin{equation}
\tilde{R}_ie^{\iota\tilde{\Psi}_i} = \frac{1}{N} \sum_j \frac{e^{\iota\theta_j}}{\omega_j-\omega_i}\, ,
\end{equation}
defines a modified order parameter with weights inversely proportional to the differences between natural frequencies. From~\cref{eq:pairwise_reduction} one can draw the following conclusions; since the control is proportional to $K_1^2$ and $R$, it tends to be large when the system has a strong coupling and it is synchronized, while it will be small when the system is desynchronized, which happens in general when $K_1\ll1$. Therefore, the control term tends to be large only when necessary. Second, the presence of $\tilde{R}_i$ induces the control injected in node (oscillator) $i$ to be stronger stronger if the natural frequency associated to this node, is close to the ones of the neighboring nodes (oscillators). Third, the cosine function allows to interpret the control signal as the original signal shifted by a quarter of period.

We can derive a similar results also for the higher-order control $h^{(M)}_i$ (see \cref{App:FullControlTerm}). The main observation is that in the present case many more terms are present, however they are proportional to $K_1^2$, $K_1K_2$ or $K_2^2$. Thus they increases if the coupling strengths are large, i.e., if the system is likely to synchronize. Because of the presence of linear integer combinations of frequencies at denominators, the control term also tends to be large if the natural frequencies are close to resonance. Finally, most of these terms  contain cosine functions, hence resulting again in shift by one fourth of period. 
Let us conclude by noticing that the proposed method could straightforwardly be adapted to more general HOKM with interaction of arbitrary order $d\geq 2$ (see Appendix~\ref{app:OrderD}).
}

\vspace{0.5cm}

{\rev 
\section{Results}
\label{sec:results}
}

The aim of this section is to perform numerical dedicated simulations to validate the control methods above introduced. More precisely, we compare the level of synchronization achieved in the HOKM~\cref{eq:HOKM3} with its controlled version given by~\cref{eq:ControlledHOKM3} and compare the results according to the used control term ($h^{(M)}$ or $\tilde{h}^{(M)}$) where we recall $M$ fixes the number of controlled nodes.
We measure the level of synchronization of the system with the standard Kuramoto order parameter {\rev $R(t)$ defined in~\cref{eq:order_parameter}}
for both the uncontrolled and the controlled HOKM. Let us recall that if $R(t)$ is close to $1$, angles are very close to each other at time $t$ and thus the system synchronizes. On the other hand, if $R(t)$ is small, {\rev the oscillators dynamics are incoherent with each other}. 
To capture the asymptotic behavior, we compute the average of the order parameter for a sufficiently long time interval, $T_{\mathrm{fin}}$, after a transient $T_{0}$, namely $\hat{R}:=\langle R(t) \rangle_{T_{0}< t< T_{0} + T_{\mathrm{fin}}}$. A good control scheme should achieve values of $\hat R$ close to zero.

{\rev In this section  will firstly assume the underlying coupling to be all-to-all for both the pairwise and the $3$-body interactions. Namely, each one of the $N$ nodes is pairwise connected to the remaining $N-1$ nodes and it participates to all possible triangles involving any distinct couples of nodes, namely, the adjacency tensors verify, for $i,j,k=1,\dots,N$,
\begin{equation}
    A_{ij}=1-\delta_{ij} \text{ and } B_{ijk}=1-\delta_{ij}\delta_{ik}\delta_{jk}
\end{equation}
where $\delta$ is the Kronecker symbol.} We will show that the {\rev full} control term $h^{(N)}$ acting on all the {\rev $N$} nodes and considering both pairwise and $3$-body interactions, is able to desynchronize the system in all the performed simulations, while $\tilde{h}^{(N)}$, i.e., still acting on all the nodes but considering only pairwise interactions, is sufficient to achieve {\rev desynchronization}, as soon as $K_1$ is not too small compared to $K_2$ {\rev (\cref{ssec:all2all})}. In a successive step we study the impact of the number of controlled nodes, $M$, on the control efficiency and show that $\hat{R}$ decreases with $M$, reaching its best performance at $M\approx 3N/5$ (\cref{ssec:number_f_controllers}). {\rev We will finally validate the proposed method by using other higher-order networks such as random simplicial complexes (\cref{ssec:Random}) and a cat connectome (\cref{ssec:cat_brain}).}

Note that triadic interactions can induce states other than full synchronization and impact basins of attractions~\cite{zhang2023deeper}. In \cref{App:Basins}, we show that for large $K_2$ values, $2$-cluster states appear but they are very unbalanced and hence close to full synchronization; let us observe that the control is effective in reducing the synchronization in this case as well.


In the following simulations, unless otherwise specified, we used higher-order structures composed of $N=50$ nodes, the natural frequencies $\omega_i$ are randomly drawn from a uniform distribution, $\omega_i\sim U([0,1])$. The initial phases, $\theta_i(0)$, are drawn from a uniform distribution close to synchronization $U([0,0.3])$. Finally, the time interval used to compute $\hat{R}$ is $[T_0,T_{\mathrm{fin}}]=[30,40]$ and we use a Runge-Kutta integrator of order 4 with a fixed integration step $0.1$ (the interested reader can refer to Appendix \ref{app:N=100} for analogous results with $N=100$). {\rev The demo codes and used hypergraphs data are available on \cite{GitLab}.}

\begin{figure*}[t]
	\centering
    \includegraphics[width=0.7\linewidth]{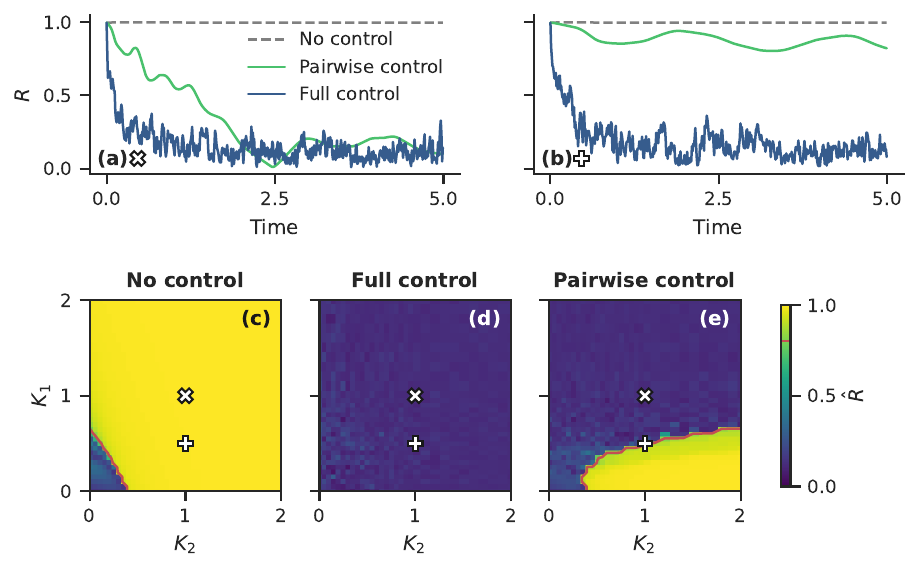}
	\caption{\rev \textbf{Desynchronizing by controlling higher-order interactions in an all-to-all hypergraph}.  We show the order parameter $R(t)$ over time for (a) $K_1=1$ and (b) $K_1=0.5$, in three cases: uncontrolled (\cref{eq:HOKM3}), with only pairwise control $\tilde{h}^{(N)}$ (green), and with full control $h^{(N)}$ (blue). Other parameters are $N=50$ and $K_2=1$.     
    The average order parameter $\hat R$ is shown over parameter space in those three cases: (c) no control, (d) full control, and (e) pairwise control. 
    The border of the synchronized region is highlighted by a level curve ($\hat{R}=0.8$) in red. The white cross and plus sign indicate the parameters used in (a) and (b), respectively. 
    }
	\label{Fig:All-to-all50}
\end{figure*}

\vspace{0.5cm}

\subsection{Desynchronizing all-to-all higher-order system}
\label{ssec:all2all}

Let us initially consider $M=N$, i.e., the control term is injected in all the nodes of the systems. Here we compare the synchronization level, measured by $R(t)$ or $\hat{R}$, as a function of the coupling strength $K_1\in [0,2]$ and $K_2\in [0,2]$ for the uncontrolled~\cref{eq:HOKM3} and the controlled systems~\cref{eq:ControlledHOKM3} by using the control terms $\tilde{h}^{(N)}$ and $h^{(N)}$.

In \cref{Fig:All-to-all50}(a-b) one can observe that the control term $h^{(N)}$ can effectively desynchronize the system. For both choices of the coupling parameters, $(K_1,K_2)=(1,1)$ and $(K_1,K_2)=(0.5,1)$, the order parameter $R(t)$ decreases rapidly for the fully controlled system (blue) while it stays close to $1$ in the uncontrolled system (gray). The pairwise control (green) manages to desynchronize the first case but not the second one.

\Cref{Fig:All-to-all50}(c-e) report the values of the averaged order parameter $\hat{R}$ in function of $K_1$ and $K_2$. We observe that the uncontrolled system exhibits strong synchronization whenever $K_2$ and $K_1$ are above some critical values (\cref{Fig:All-to-all50}(c)), which is consistent with previous results stating that higher-order interactions strengthen the local stability of the synchronized state~\cite{skardal2019abrupt,Carletti2020dynamical,gambuzza2021stability,millan2020explosive,zhang2023deeper}. The red curve indicates the parameters for which $\hat{R}=0.8$, i.e., an arbitrary (large enough) value that we fixed to define a sufficiently large level of synchronization; only for small enough $K_1$ and $K_2$, the system does not synchronize (see blue region in the bottom-left corner of panel c). Let us observe that the region where $\hat{R}>0.8$ covers {\rev approximately $96\%$} of the considered domain.

One can furthermore remark that the second-order interaction strength $K_2$ seems to play a more important role in synchronization than $K_1$. Indeed, in \cref{Fig:All-to-all50}(c) for $K_1=0$ the value {\rev $K_2\approx 0.39$} is sufficient to synchronize the uncontrolled system, indeed passing this value $\hat{R}$ suddenly passes from very low values to very large ones. On the other hand for $K_2=0$ the system requires {\rev $K_1\approx 0.66$} to synchronize. 

On the other hand, the fully controlled system (see \cref{eq:ControlledHOKM3}) by using $h^{(N)}$ as control term~\cref{eq:h_i}, is successfully prevented from synchronizing for all $K_1$ and $K_2$ in the considered range (\cref{Fig:All-to-all50}(d)). The mean value of $\hat{R}$ over the whole domain is $0.197$ 
and it never reaches the threshold $0.8$ (its maximal value is approximately 0.25).  This indicates that the proposed control scheme is successful in reducing synchronization.

{\rev In \cref{Fig:All-to-all50}(e),} we consider the control term $\tilde{h}^{(N)}$ (see \cref{eq:h_iPartial2D}) in which only the pairwise part is taken into account. 
We can appreciate that the control method is remarkably efficient to reduce synchronization being the values of $\hat{R}$ as low as the ones obtained in Fig.~\ref{Fig:All-to-all50}(d) for most of the $(K_1,K_2)$ choices, especially if $K_1>K_2$. As $\tilde{h}^{(N)}$ is less complicated to compute and involves a significantly lower cost in terms of total energy of the injected signal (see Appendix \ref{app:Cost} for more details) the option of using this lighter version in some context could be very interesting for future applications. 

Nevertheless the pairwise control $\tilde{h}^{(N)}$ is not sufficient to desynchronize the system if $K_1$ is small in comparison to $K_2$. Indeed we can roughly observe than if $K_1<0.5$ and {\rev $K_2>0.39$} then $\hat{R}>0.8$ showing that in this case the control with only pairwise terms is not able to desynchronize the system, indeed the yellow region then covers {\rev approximately $21\%$} of the considered range of parameters. These latter observations are consistent with \cref{Fig:All-to-all50}(a-b) where we show that $R(t)$ barely decreases by using $\tilde{h}^{(N)}$ (green) with $(K_1,K_2)=(0.5,1)$ while it goes below $0.2$ if $(K_1,K_2)=(1,1)$.

\begin{figure}[t!]
    \centering
    \includegraphics[width=0.8\linewidth]{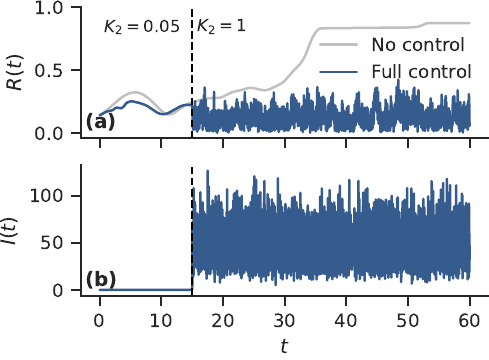}
    \caption{{\rev \textbf{The control is adaptive}. We show (a) the evolution of the order parameter, $R(t)$, in the system with no control (gray) and with full control $h^{(N)}$ (blue), and (b) the intensity of the full control over time.  
    Initially, we set both coupling strengths to small enough values, $K_1=K_2=0.05$: the systems---both uncontrolled and controlled---do not synchronize $R(t)\approx 0$ in the whole time interval, and the intensity of the control is null $I(t)=0$ (b). At $t=15$, we set $K_2=1$ so that the uncontrolled system synchronizes. The  control intensity increases and prevents the controlled system from synchronizing. }
    The initial angles are drawn from an uniform distribution, $\theta(0)\sim U([0,2\pi])$.}
    \label{Fig:K2Bifurcation}
\end{figure}

{\rev Before concluding this part, we briefly analyze the control magnitude in function of the systems parameters.} As already observed in \cref{sec:control}, the control intensity is proportional to the coupling strengths $K_1$, $K_2$ and inversely proportional to the frequency differences. In particular, the control intensity will be higher when $K_1$ and $K_2$ are large (see also \cref{app:Cost}), setting for which the local synchronization basin is larger. There is therefore a predominant dynamical effect: the control is large whenever the uncontrolled system would tend to synchronization. This fact is illustrated in \cref{Fig:K2Bifurcation} {\rev where we report the time evolution of the order parameters $R(t)$}; once we set the coupling strengths to $K_1=K_2=0.05$, both the uncontrolled and the controlled systems will evolve toward an asynchronous regime (see panel (a)). {\rev Indeed $R(t)$ results to be quite small, meaning that the system} is far from synchronization. At the same time, the control intensity, measured by {\rev the nodes average of the absolute value of the control signal}, $I(t):=\langle |h^{(N)}_i| \rangle_{i=1,\dots,N}$ (panel (b)), remains close to zero for a large time interval. Then at $t=15$ we modify the coupling parameter $K_2$ and we set $K_2=1$, the uncontrolled system (gray line in panel (a)) will shortly after synchronize, indeed the order parameter will steadily increase while the controlled system will remain in an asynchronous state (blue line in panel (a))). To achieve this, the controlled system requires a strong enough control intensity as one can appreciate by looking at  the bottom panel.

{\rev In conclusion,} we can thus claim that the control term $h^{(N)}$ is able to desynchronize the HOKM even if there are only higher-order interactions involved and thus generalizing the results presented in~\cite{gjata2017using}. Moreover, the proposed control method is robust to the superposition of interaction of different orders, despite the fact that the stability of the synchronization state is strengthen once higher-order terms are taken into account~\cite{zhang2023deeper}. {\rev Furthermore,} in many cases, we can prevent a higher-order system from synchronizing by using only the control obtained by solely using the pairwise interaction, {\rev $\tilde{h}_i^{(N)}$}, resulting in a huge computational advantage; however if the pairwise coupling is weak enough with respect to the second-order one, then the full higher-order control---with pairwise and triadic terms---is needed to achieve desynchronization.

\subsection{Impact of the number of controllers}
\label{ssec:number_f_controllers}

The aim of this section is to study the impact of the number of controllers $M$ on the control efficiency, i.e., how $\hat{R}$ depends on $M$, again by using both control terms $h^{(M)}$ and $\tilde{h}^{(M)}$. Here again we restrict our analysis to the case of all-to-all pairwise and $3$-body interactions.

\begin{figure}[t]
	\centering
	\includegraphics[width=\linewidth]{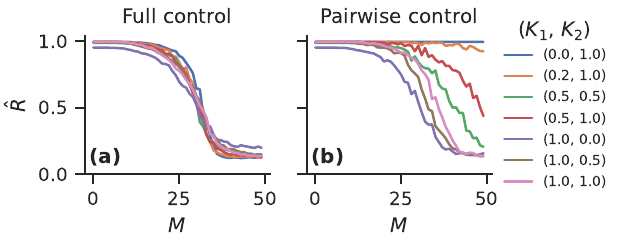} 
	\caption{{\rev \textbf{Impact of the number of controlled nodes.} We report the average Kuramoto order parameter, $\hat{R}$, as a function of the number $M$ of controlled nodes for several values of $K_1$ and $K_2$, with (a) full control $h^{(M)}$, and (b) pairwise control $\tilde{h}^{(M)}$.
    Each curve represents the average of $150$ independent numerical simulations with different samples for $\pmb{\omega}\sim U([0,1])$ in an all-to-all hypergraph with $N=50$ nodes.}
    }
	\label{Fig:All-to-all50:Partial}
\end{figure}

Let us observe that in this setting all nodes, i.e., oscillators, are equivalent regarding the higher-order network topology and thus they differ only for the natural frequencies. We can thus safely assume that the number $M$ of controllers is the key parameter and not their position in the structure. This would not longer be true for general hypergraphs topologies where two different subsets of controllers of the same size $M$, may not return similar outcomes because of the position of the controllers nodes in the network, as already shown in the study of complex networks~\cite{wang2002pinning,Wenwu2009onpinning,liu2011controllability,liu2016control,liu2021optimizing,Moriame2023onthelocation} and, more recently, on higher-order networks~\cite{delellis2023pinning,wang2024synchronization,Muolo2024Pinning,Li2024synchronization,Chen2021controllability}. In those cases, the control efficiency depends on particular characteristics of the (higher-order) networks topologies and the centrality scores of the selected nodes (very often their degree). The identification of the optimal pinned subset is thus a whole problem in itself and we will not consider it in this work.

\begin{figure*}[t]
    \centering
    \includegraphics[width=0.8\linewidth]{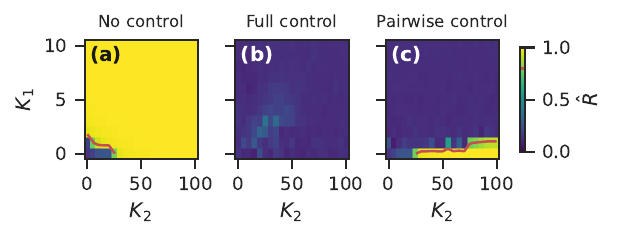}
    \caption{\rev \textbf{The case of random simplicial complexes.} We show the average order parameter $\hat{R}$ over parameter space $(K_1,K_2)\in[0,10]\times[0,100]$ in three cases: 
    (a) no control, (b) full control  $h^{(N)}$, and (c) pairwise control $\tilde{h}^{(N)}$. 
    The random simplicial complex has $50$ nodes and average degrees $\left< k_1\right>=40$ and $\left<k_2\right>=20$.
    The border of the synchronized region is highlighted by a level curve ($\hat{R}=0.8$) in red.  }
    \label{Fig:RandomTopologies}
\end{figure*}

In \cref{Fig:All-to-all50:Partial} we show $\hat{R}$ as a function of the number of controlled nodes $M$ for several values of $K_1$ and $K_2$. The first observation is that $\hat{R}$ significantly decreases with increasing $M$ once we use the control term built by using both pairwise and $3$-body interactions, independently from the choices of $(K_1,K_2)$ (see \cref{Fig:All-to-all50:Partial}(a)). Namely, the larger the number of pinned nodes the better is the control to achieve desynchronization. Moreover the different functional forms of $\hat{R}$ versus $M$ are very similar each other, even for different values of $K_1$ and $K_2$. $\hat{R}$ remains constant at its maximal value until $M\approx N/5$ then it rapidly decreases before stabilizing to a lower value. The latter is obtained for $M$ between $3N/5$ and $4N/5$ in the case of the full control $h^{(N)}$. 

In \cref{Fig:All-to-all50:Partial}(b) we present the results obtained by using the control built by using only pairwise terms. In many considered cases, the decreasing behavior of $\hat{R}$ versus $M$ is the same of the one presented in~\cref{Fig:All-to-all50:Partial}(a). There are however notable exceptions for $K_2=1$ and $K_1\in [0,0.5]$, indeed when $K_1=0.2$, the control is still not able to desynchronize the system, $\hat{R}$ stays larger than $0.9$ for all $M$ and barely does not decrease. Finally when $K_1=0.5$, $\hat{R}$ decreases with $M$ but much slower than in other cases and cannot reach a value of $0.4$ or smaller. In conclusion, for those values of $K_1$ and $K_2$, the control term $\tilde{h}^{(M)}$ cannot desynchronize the system. Note also that for $K_1=K_2=0.5$ the decreasing is fast but still cannot reach the same lower values than in  \cref{Fig:All-to-all50:Partial}(a). Also in this case $\tilde{h}^{(M)}$ could not be considered sufficient to reduce synchronization.

{\rev
\subsection{Desynchronization in random simplicial complexes}
\label{ssec:Random}
}

\begin{figure*}[hbt]
    \centering
    \includegraphics[width=0.8\linewidth]{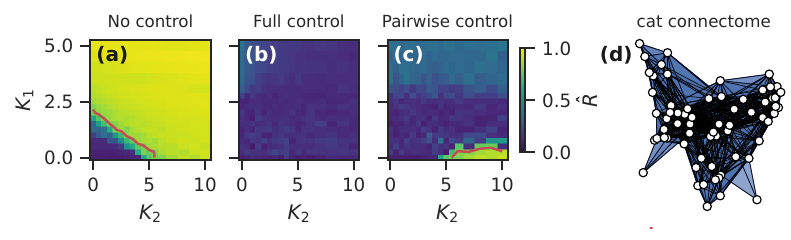}
    \caption{\rev \textbf{Desynchronizing a brain connectome.} We show the average order parameter $\hat{R}$ over parameter space $(K_1,K_2)\in[0,5]\times[0,10]$ in three cases: 
    (a) the original model, i.e., without any control, (b) the full control  $h^{(N)}$, and (c) the pairwise control $\tilde{h}^{(N)}$. 
    The boundary of the synchronized region is highlighted by the level curve ($\hat{R}=0.8$) in red. 
    (d) Visual representation of the hypergraph. The higher-order brain connectome has $N=65$ nodes, $730$ links and $3613$ triangles.}
    \label{fig:cat_brain}
\end{figure*}

We now demonstrate {\rev the efficiency of the proposed control method also} on less regular structures, {\rev such as,} random simplicial complexes, built following the algorithm proposed in~\cite{iacopini2019simplicial}.
The {\rev latter creates} a $2$-simplicial complex with an \textit{a priori} fixed mean node degree, $k_1$, and mean hyperdegree, $k_2$, where the hypedregree of node $i$ denotes the number of triangles incident with node $i$. 

{\rev Let us briefly describe the main structure of the algorithm and invite the interested reader to consult the original work~\cite{iacopini2019simplicial}. First, for every triplet of distinct nodes $(i,j,k)$ we add with probability $\frac{2k_2}{(N-1)(N-2)}$ a triangle formed by the latter nodes so that the mean hyperdegree $k_2$ is reached in expectation. Then, the three links formed by the couples, $(i,j)$, $(j,k)$ and $(k,i)$, are added. At this point the average node degree is not necessarily equal to $k_1$ in expectation, to achieve also this goal, supplementary links are added connecting pairs of distinct nodes $(i,j)$ with probability $\frac{k_1-2k_2}{N-1-2k_2}$.}

In Fig.~\ref{Fig:RandomTopologies} we show the average order parameter $\hat{R}$ for the uncontrolled system \cref{eq:HOKM3} coupled by using a random $2$-simplicial complex generated according to method stated above with $N=50$ nodes, average degree $k_1=40$ and average hyperdegree $k_2=20$, and its controlled version \cref{eq:ControlledHOKM3} by using $M=N$, i.e., all are nodes controlled, in the version with only pairwise terms and the one allowing for both pairwise and $3$-body interactions.

First of all, we can observe that the results agree with those previously presented, i.e., the uncontrolled system synchronizes when $K_1$ or $K_2$ is sufficiently large, and the system controlled by $h^{(N)}$ completely removes the synchronization. Finally, the control restrained to pairwise interactions is also able to desynchronize the system provided $K_1$ is not too small and for (almost) any $K_2$. Note that here the used values for $K_1$ and $K_2$ are larger than then ones previously used (see Fig.~\ref{Fig:All-to-all50}). This is due to the lower density of {\rev the interactions in those } simplicial complexes{\rev, as not all links and triangle are present, and the normalizations} by $N$ and $N^2$. One thus needs larger interaction strength so that the trajectories stay close to the synchronization. The region where the pairwise control is not capable to reduce synchronization is about $4.1\%$ of the whole parameters set, this corresponds to the area where $\hat{R}>0.8$ (in yellow in the figure), and it is smaller than the similar region in the case of all-to-all (see Fig.~\ref{Fig:All-to-all50}).

\vspace{0.5cm}

{\rev
\subsection{Desynchronization in a cat brain connectome}
\label{ssec:cat_brain}

As another example of the applicability of the method on a non-regular higher-order topology, we use as underlying coupling scheme a cat brain connectome~\cite{networkswebsite}. Let us however observe that the latter was available only with pairwise connections, i.e., encoded by the symmetric adjacency matrix $\mathbf{A}$, we thus decided to promote to $3$-body interaction all the available triangles, namely for all triplets $(i,j,k)$ such that $A_{ij}=A_{jk}=A_{ki}=1$ we set $B_{\pi(ijk)}=1$, where again $\pi$ denotes all the permutations of the three indexes. The resulting simplicial complex has $N=65$ nodes, $730$ links and $3613$ triangles ($3$-body hyperedges). 

The results are reported in Fig.~\ref{fig:cat_brain} and are in line with those shown in Figs.~\ref{Fig:All-to-all50} and~\ref{Fig:RandomTopologies}. Indeed, the full control term $h^{(N)}$ desynchronizes the system for every tested couples $(K_1,K_2)\in[0,5]\times[0,10]$ and so does $\tilde{h}^{(N)}$ except for values in a (relatively) small region corresponding to small $K_1$ ($K_1<0.3$) and large $K_2$ ($K_2\leq 5.5$) (see the bottom right corner) for which the pairwise control is not capable to desynchronize the system.

Let us conclude this analysis with the following observation; the efficiency of the pairwise control $\tilde{h}^{(N)}$ slightly decreases for large values of $K_1$ ($K_1>2.5$) independently from $K_2$, indeed $\hat{R}$ slightly increases. This effect is not observed in the case of the full control.}

\section{Discussions and conclusions}
\label{conclusion}

In this work, we proposed a generalization of the feedback pinning control method developed in~\cite{gjata2017using,asllani2018minimally} in order to desynchronize the higher-order Kuramoto model where $3$-body interactions coexist with pairwise ones. The strategy relies on the implementation of a pinning control term obtained by applying Hamiltonian control theory~\cite{vittot2004perturbation,ciraolo2004control} to a suitable embedding Hamiltonian system.

The Hamiltonian function we propose is a natural generalization of the one used in~\cite{witthaut2014kuramoto}, to deal with higher-order interactions. Let us first remark that the embedded HOKM contains a combination of the two possible types of second-order interactions \cref{3orderTerm1} and \cref{3orderTerm2} often proposed separately in the literature dealing with extensions of the Kuramoto pairwise model. They appear together in the context of phase reduction studies of coupled oscillators~\cite{leon2019phase,gengel2020highorder}; the derivation of the Hamiltonian system we propose is another, independent, case in which this combination of second-order interactions terms arise naturally. Moreover, although in most of the available literature they are considered separately, there is no clear consensus on which one should be used in function of the context and neither any clear argument that would justify that they even have to be separated. Hence, we believe that the combined system \eqref{eq:HOKM3} would deserve deeper investigations in further works.

We numerically showed that the proposed control term acting on all the nodes and taking into account both pairwise and $3$-body interactions, enables to desynchronize the HOKM system~\eqref{eq:HOKM3} {\rev on structures of increasing complexity: all-to-all higher-order networks, random simplicial complexes, and a higher-order cat connectome}. The proposed theory supports the claim that this result remains valid on a wider set of higher-order networks. For all considered values of interaction strengths $K_1$ and $K_2$ the smallness of averaged order parameter $\hat{R}$ proves that the system is desynchronized while the uncontrolled one is, on the contrary, strongly synchronized. We also numerically showed that the control intensity is larger when the parameters are favorable to synchronization in the uncontrolled system which makes the control efficient and less invasive in any situation.

We also investigated the possibility of using a control term based on the pairwise interaction terms only, as done in~\cite{gjata2017using}, to desynchronize the HOKM. Interestingly, the strategy works if the pairwise coupling strength, $K_1$, is not too small  {\rev while} the higher-order coupling strength, $K_2$, {\rev is large}. In this case, the perturbation induced by the pairwise part {\rev $\tilde{h}^{(M)}$} of the control term is sufficient to desynchronize the system and has a lower cost
than $h^{(M)}$. On the other hand, when $K_1$ {\rev lies below a certain threshold, itself depending on $K_2$}, it is necessary to consider the second-order interaction term in the control, and thus the higher-order generalization is mandatory. {\rev The above mentioned threshold seems also to strongly depend on the underlying topology, as we observed for the all-to-all hyperpraghs, random and brain simplicial complexes.}

Moreover we investigated the efficiency of the control $h^{(M)}$ as a function of the number of pinned nodes $M$. It is clear that $\hat{R}$ decreases with increasing $M$, for all values of $K_1$ and $K_2$. The functional form of $\hat{R}$ versus $M$ is quite robust with respect to the used values of $K_1$ and $K_2$, it shows a horizontal plateau for small $M$, before decreasing quite rapidly and then reaching a second plateau at its minimal value for large $M$. 

Let us observe that the proportion of controlled nodes required to reach the minimal value of $\hat{R}$ is approximately $0.6$ that results to be a quite large value in a pinning control perspective. This can be explained by noticing that the stability of the synchronized state is intrinsically very strong in the all-to-all topology. One therefore needs a high intensity of control to move trajectories far away from the synchronization state. In more sparse systems, the required proportion of controlled nodes to get desynchronization should be lower; moreover for heterogeneous higher-order structures, the determination of the optimal pinning nodes is an open problem, still in its infancy. For this reason we believe this analysis goes beyond the scope of this work and will be investigated elsewhere. 

Another perspective of this work is to, somehow, simplify the control function in~\cref{eq:h_i}, that contains many terms and requires a deep understanding of the system, {\rev as it has been done in~\cite{asllani2018minimally} for the pairwise system. The advantages of this analysis is twofold. It could help to optimize the control by keeping} only the ``more'' effective terms. Such simplified form could be relevant for possible practical applications, e.g., epilepsy and neural diseases, or in the general theory of control of higher-order interactions. {\rev It could also help to get a deeper understanding of the control term, which remains complex and whose dynamics are hard to understand. 

Furthermore, further study should carefully investigate the threshold for $K_1$ that guarantees the desynchronization by using the pairwise term alone. As already observed, the latter should depend on $K_2$ and  the underlying topology.

Finally, a remaining interesting research line is the investigation on a general method to build an embedding Hamiltonian systems, when possible. Indeed, an essential assumption for the use of Hamiltonian control theory on a synchronizing system is either that the latter is Hamiltonian, which is in general not the case, or that it exists an embedding system which is Hamiltonian. This is the case for (higher-order) Kuramoto model. So it is for the Stuart-Landeau model 
\begin{equation}\label{start_landau}
    \dot{z_i}= \sigma z_i -\beta z_i|z_i|^2 + \frac{K_1}{N}\sum_j A_{ij}z_j\, ,
\end{equation}
whose complex oscillators can exhibit phase and/or amplitude oscillations, with the complex Hamiltonian function $H_{SL}(z,\bar{z}) = \sum_{i}\imath(\sigma|z_i|^2 -\frac{\beta}{2}|z_i|^4) +\imath \frac{K_1}{N}\sum_{i,j}A_{ij}\frac{|z_j|^2}{2}$. This can be used as a basis for the development of a desynchronizing control method for the (higher-order) Stuart-landau model However, to the best of our knowledge, there is no systematic method to construct such Hamiltonian function. If such method was found, the Hamiltonian control theory could eventually be applied to more diverse models of oscillators.

}

\section*{Acknowledgements}

Part of the results were obtained using the computational resources provided by the ``Consortium des Equipements de Calcul Intensif" (CECI), funded by the Fonds de la Recherche Scientifique de Belgique (FRS-FNRS) under Grant No. 2.5020.11 and by the Walloon Region. M.L. is a Postdoctoral Researcher of the Fonds de la Recherche Scientifique–FNRS with project Under-Net 40016866.

\bibliography{bib}

\appendix
\onecolumngrid

{\rev

\section{Derivation of the control term}

The aim of this section is to provide more details about the derivation of the control term used in the main text.

\subsection{About the Hamiltonian system}
\label{app:Embeding}

We hereby prove that the higher-order Hamiltonian system~\eqref{def:HOHam3} allows to recover the higher-order Kuramoto model \cref{eq:HOKM3} once restricted on the family of invariant tori $T_c:=\{(\mathbf{I},\pmb{\theta})\in\mathbb{R}_+^N\times[0,2\pi]^N | \forall i\in\{1,\dots,N\}:I_i=c\}$.

Let us thus derive the action and angle dynamics. First, we have
\begin{eqnarray}
	\dot{I}_i&=&-\frac{\partial H}{\partial \theta_i} \\
    &=&-2\frac{K_1}{N}\sum_j A_{ij}\sqrt{I_iI_j}(I_j-I_i)\cos(\theta_j-\theta_i)+ \nonumber \\
	&&-2\frac{K_2}{N^2}\sum_{j,k}B_{ijk}\sqrt[3]{I_iI_jI_k}(I_j+I_k-2I_i)\cos(\theta_j+\theta_k-2\theta_i)+ \nonumber \\
	&&-2\frac{K_2}{N^2}\sum_{j,k}B_{ijk}\sqrt[3]{I_iI_jI_k}(2I_j-I_k-I_i)\cos(2\theta_j-\theta_k-\theta_i)\, , \nonumber 
\end{eqnarray}
from which it is clear that the torus $T_c$ is invariant for all $c>0$, indeed by inserting $I_i=c$ into the latter we get $\dot{I}_i=0$ and thus the action variables will not evolve. This property will play a relevant role in the following, as it was the case for the results presented in~\cite{witthaut2014kuramoto}.

Then the angles evolution is given by
\begin{eqnarray}
	\dot{\theta_i}&=&\frac{\partial H}{\partial I_i} \\
    &=&\omega_i- 2\frac{K_1}{N}\sum_j A_{ij}\sin(\theta_j-\theta_i)\left[\frac{1}{2}\sqrt{\frac{I_j}{I_i}}(I_j-I_i)-\sqrt{I_jI_i}\right]+ \nonumber \\
	&&-\frac{K_2}{N^2}\sum_{j,k}B_{ijk}\sin(\theta_j+\theta_k-2\theta_i)\left[-2\sqrt[3]{I_iI_jI_k} +\frac{1}{3}\sqrt[3]{\frac{I_jI_k}{I_i^2}}(I_j+I_k-2I_i) \right]+ \nonumber \\
	&&+2\frac{K_2}{N^2}\sum_{j,k}B_{ijk}\sin(2\theta_j-\theta_k-\theta_i)\left[\sqrt[3]{I_iI_jI_k} -\frac{1}{3}\sqrt[3]{\frac{I_jI_k}{I_i^2}}(2I_j-I_k-I_i) \right] \, , \nonumber 
\end{eqnarray}
and when we set $I_i=\frac{1}{2}$ for all $i\in\{1,\dots,N\}$, we finally obtain~\cref{eq:HOKM3}. Observe here again that by taking $I_i=c\neq\frac{1}{2}$ will also reduce to~\cref{eq:HOKM3} but with rescaled interaction strengths $2cK_1$ and $2cK_2$.

Moreover, we can see that the synchronized regimes of the HOKM \cref{eq:HOKM3} correspond to instable behaviors of the invariant Torus $T_{1/2}$ for the  Hamiltonian system. Indeed let us compute the Jacobian matrix $\mathbf{J}(\mathbf{I},\pmb{\theta})$ of the system evaluated on the torus, $\mathbf{I}_i=\frac{1}{2}$, and we assume the angles to be very close to synchronization, namely they are very close to $\mathcal{S}_{\tilde{\theta}}=\{ \theta_i \in [0,2\pi]: \theta_i-\tilde{\theta}=0\}$, for a generic $\tilde{\theta}$. A straightforward computation allows us to obtain
\begin{equation}\label{eq:Jacobian}
	\mathbf{J}=
	\begin{bmatrix}
		\mathbf{L}&0\\0&-\mathbf{L}
	\end{bmatrix}
\end{equation}
where the $N\times N$ matrix $\mathbf{L}:=2\mathbf{L}^{(1)}+6\mathbf{L}^{(2)}$ is given by
\begin{equation}\label{def:multiorderlaplacian}
	\begin{split}
		L^{(1)}_{ij}=\left\{
		\begin{matrix}
			-\frac{K_1}{N}A_{ij}&\text{  if }i\neq j\\
			-\sum_{j=1}^N L^{(1)}_{ij}&\text{  if }i=j
		\end{matrix}
		\right.\, ~~\text{and}
		&~~
		L^{(2)}_{ij}=\left\{
		\begin{matrix}
			-\frac{K_2}{N^2}\sum_{k=1}^N B_{ijk}&\text{  if }i\neq j,\\
			-\sum_{j=1}^N L^{(2)}_{ij}&\text{  if }i=j.
		\end{matrix}
		\right.
	\end{split}
\end{equation}
We note that $\mathbf{L}$ has the same form as the multiorder Laplace matrix introduced in~\cite{lucas2020multiorder} as we perform the sum on the third index of $\mathbf{B}$ to obtain a two-dimensional matrix to be added to $\mathbf{A}$.

The matrix $\mathbf{L}$ is by construction symmetric, non-negative and therefore has eigenvalues $0=\lambda_1<\lambda_2\leq\dots\leq\lambda_N$. In conclusion $\mathbf{J}$ has two null eigenvalues, that characterize the fact that $T_c$ and $\mathcal{S}_{\tilde{\theta}}$ are (almost) invariant. The negative eigenvalues $-\lambda_2,\dots,-\lambda_N$ associated to the direction tangent to the torus, force the angle variables to stay close to $\mathcal{S}_{\tilde{\theta}}$; the positive eigenvalues associated to the orthogonal direction to the torus, on the other hand tend to move the orbit far away from the torus. In other words, when HOKM is synchronized then $T_c$ is an unstable invariant manifolds for~\eqref{def:HOHam3} for $c>0$.
}

{\rev

\subsection{About the Hamiltonian control theory}
\label{App:GeneralControlTerm}

In this Section we briefly describe the main steps required to obtain the control $f(V)$ and we ref the interested readers to consult~\cite{vittot2004perturbation,ciraolo2004control,gjata2017using}. We here build a perturbation $f(V)$ of the Hamiltonian function $H$ defined in~\eqref{def:HOHam3}. The latter should be chosen in such a way the Hamiltonian flow induced by $H + f(V)$ is canonically conjugate to the one induced by $H_0 + G(V)$ where the function $G(V)$ is the resonant part of $V$. Under the non-resonance condition on the frequencies hypothesized in the main text, we get $G(V)=0$.

Let $\mathcal{A}$ be the Lie algebra formed by $C^{\infty}$ functions of action-angle variables $(\mathbf{I},\pmb{\theta})$ in $\mathbb{R}_+^N\times[0,2\pi]^N$ with values in $\mathbb{R}^{2N}$ and by the Poisson brackets $\{\cdot,\cdot\}$ defined by
\begin{equation}
	\forall f,g\in\mathcal{A}: \{f,g\}:=\frac{\partial f}{\partial \mathbf{I}}\cdot\frac{\partial g}{\partial\pmb{\theta}} - \frac{\partial f}{\partial \pmb{\theta}}\cdot\frac{\partial g}{\partial\mathbf{I}}\, ,
\end{equation}
where $\cdot$ denotes the (real) scalar product. For any $f\in\mathcal{A}$, we can define the linear operator
\begin{eqnarray}
\label{eq:opf}
 &\{f\}&:\mathcal{A}\longrightarrow\mathcal{A}\notag\\
 &g&\mapsto\{f\}g:=\{f,g\}\, .
\end{eqnarray}

Because of the Hamiltonian structure of the system, the time evolution of any function $g\in\mathcal{A}$ is given by 
\begin{equation}
\label{eq:gevolv}
g(x)=e^{t\{H\}}g(x_0):=\sum_{n=0}^{+\infty}\frac{t^n\{H\}^n}{n!}g(x_0)\, ,
\end{equation}
where $x$ is the position of the orbit at time $t$ starting from $x_0$ at $t=0$.

Let now $\Gamma:\mathcal{A}\longrightarrow\mathcal{A}$ be the pseudoinverse operator of $\{H_0\}$, namely
\begin{equation}\label{eq:Gamma}
	\{H_0\}^2\Gamma=\{H_0\}\, .
\end{equation}
To compute one of the possible solutions of~\cref{eq:Gamma}, authors of~\cite{vittot2004perturbation,ciraolo2004control} first 
wrote $V$ by using its Fourier series, $V=\sum_{\mathbf{k}\in\mathbb{Z}^N}V_{\mathbf{k}}(\mathbf{I})e^{\iota \mathbf{k}\cdot\pmb{\theta}}$, where $\iota = \sqrt{-1}$, and then computed
\begin{equation}
\{H_0\}V=\sum_{\mathbf{k}\in\mathbb{Z}^N} \iota \, (\pmb{\omega}\cdot\mathbf{k}) V_{\mathbf{k}}(\mathbf{I})e^{\iota \mathbf{k}\cdot\pmb{\theta}}\, , 
\end{equation}
where use has been made of the definition $\frac{\partial H_0}{\partial I_i}=\omega_i$. Then one can deduce a formal  definition of $\Gamma$ through
\begin{equation} \label{eq:GammaV}
	\Gamma V= \sum_{\mathbf{k}\in\mathbb{Z}^N: \, \pmb{\omega} \cdot \mathbf{k} \neq 0} \frac{V_{\mathbf{k}}(\mathbf{I})e^{\iota \mathbf{k}\cdot\pmb{\theta}}}{\iota \,(\pmb{\omega} \cdot \mathbf{k})}
\end{equation}
which clearly satisfies
\begin{equation}
	\{H_0\}\Gamma V =V\, .
\end{equation}

Finally it has been proved~\cite{vittot2004perturbation,ciraolo2004control} that the desired perturbation can be defined as 
\begin{equation}
	f(V):=\sum_{n=1}^{\infty}\frac{(-1)^n\{\Gamma V\}^n}{(n+1)!}(n\mathcal{R}+1)V\, ,
\end{equation}
where the exponent $n$ denotes the $n$-th composition of the operator $\{\Gamma V\}$ and $\mathcal{R}$ is the resonant operator. The latter is defined as $\mathcal{R}:=\mathbf{1}-\{H_0\}\Gamma$ where $\mathbf{1}$ is the identity operator. With this definition of $f(V)$, the goal is achieved with the function $G$ defined by
\begin{equation}
G:= \mathcal{R}V=\sum_{\mathbf{k}\in\mathbb{Z}^N: \, \pmb{\omega} \cdot \mathbf{k} = 0}V_{\mathbf{k}}(\mathbf{I})e^{\iota \mathbf{k}\cdot\pmb{\theta}}\, .
\end{equation}

In the case under scrutiny one can realize from~\eqref{def:HOHam3} that the non-zero Fourier coefficients of $V$ are given by $\mathbf{k}^{(1)}_{ij}:=\mathbf{e}_i-\mathbf{e}_j$ for the pairwise term, and $\mathbf{k}^{(2)}_{ijk}:=\mathbf{e}_j+\mathbf{e}_k-2\mathbf{e}_i$ for the three-body interaction, where $\mathbf{e}_i$ is the $i^{th}$ canonical basis vector of $\mathbb{R}^N$. We can eventually rewrite $V=V^{(1)}+V^{(2)}$ and $\Gamma V=\Gamma V^{(1)}+\Gamma V^{(2)}$ by separating the two types of terms, we observe that there are no resonant terms, and thus $G=0$.

In conclusion the function 
\begin{equation}
	f(V):=\sum_{n=1}^{\infty}\frac{(-1)^n\{\Gamma V\}^n}{(n+1)!}V\, ,
\end{equation}
is a suitable control term for the system~\eqref{def:HOHam3}. As we eventually want to control its restriction on $T_{\frac{1}{2}}$, the control term $h_i$ added to the evolution equation~\cref{eq:HOKM3} for the $i$-{th} oscillator $\theta_i$ is defined by $h_i:=\frac{\partial f(V)}{\partial I_i}\big|_{\mathbf{I}=\frac{1}{2}}$.

To avoid dealing with the infinite series that defines the control $f(V)$ and thus to the problem of its convergence, we truncate the series and only keep the dominant term, i.e., the first one $n=1$ as done in~\cite{gjata2017using} to recover~\eqref{eq:h_i}.
}

\subsection{Derivation of the control term with pairwise and 3-body interactions}
\label{App:FullControlTerm}

The aim of this section is to provide the reader the details about the construction of the control terms given in \cref{eq:h_iPartial,eq:h_iPartial2D}. Let be $M\leq N$ and $\{1,\dots,M\}$ the subset of observed-controlled nodes after a proper nodes labeling (\cref{eq:h_i} will be recovered when $N=M$). 

{\rev
From Eqs.~\eqref{eq:h_i},~\eqref{eq:h_iPartial} and~\eqref{eq:split_V} we get
\begin{equation}
    \begin{split}
        h_i^{(M)}=&-\frac{1}{2}\frac{\partial }{\partial I_i}\left[\{\Gamma V^{(M,1)}\}V^{(M,1)}+\{\Gamma V^{(M,1)}\}V^{(M,2)}\right. \\
        &\left.+\{\Gamma V^{(M,2)}\}V^{(M,1)}+\{\Gamma V^{(M,2)}\}V^{(M,2)}\right]\big|_{\mathbf{I}=\frac{1}{2}}\, .
    \end{split}
\end{equation}
}
The non-zero coefficients of the Fourier series of $V^{(M)}=V^{(M,1)}+V^{(M,2)}$ are respectively given by
\begin{equation}
		V_{\mathbf{k}}^{(M,1)}=
		\left\{
		\begin{matrix}
			-\frac{K_1}{N}A_{ij}f_1(I_i,I_j)&~\text{if } \mathbf{k}=\mathbf{e}_j-\mathbf{e}_i;~i\neq j\in \{1,\dots,M\},\\
			0&~\text{otherwise}
		\end{matrix}
		\right.
\end{equation}
and
\begin{equation}
		V_{\mathbf{k}}^{(M,2)}=
		\left\{
		\begin{matrix}
			-\frac{K_2}{N^2}B_{ijk}f_2(I_i,I_j,I_k)&~\text{if } \mathbf{k}=\mathbf{e}_j+\mathbf{e}_k-2\mathbf{e}_i;~~i,j,k\in \{1,\dots,M\};i\neq j,j\neq k, k \neq i,\\
			0&~\text{otherwise}
		\end{matrix}
	\right.
\end{equation}
where $f_1(I_i,I_j):=\sqrt{I_iI_j}(I_j-I_i)$ and $f_2(I_i,I_j,I_k):=\sqrt[3]{I_iI_jI_k}(I_j+I_k-2I_i)$. By using the definition of the pseudo-inverse operator $\Gamma$ and by injecting the latter functions in \cref{eq:GammaV} we get
\begin{eqnarray}\displaystyle
		\Gamma V^{(M,1)}&= &\sum_{i,j=1}^MA_{ij}\frac{f_1(I_i,I_j)}{\omega_j-\omega_i}\cos(\theta_j-\theta_i),\\ 
		\Gamma V^{(M,2)}&= &\sum_{i,j,k=1}^MB_{ijk}\frac{f_2(I_i,I_j,I_k)}{\omega_j+\omega_k-2\omega_i}\cos(\theta_j+\theta_k-2\theta_i)\, .
\end{eqnarray}

Let us observe that the functions $f_1$ and $f_2$ are the only part of $V$ and $\Gamma V$ that depend on $\mathbf{I}$. If all the action variables are set to the same constant value $c=\frac{1}{2}$ then $f_1$ and $f_2$ vanish and their derivatives are given by 
{\rev
\begin{eqnarray}
    \frac{\partial f_1(I_i,I_j)}{\partial I_j}\rvert_{I_i=1/2}=-\frac{\partial f_1(I_i,I_j)}{\partial I_i}\rvert_{I_i=1/2}&=&1\\
    \frac{\partial f_2(I_i,I_j,I_k)}{\partial I_j}\rvert_{I_i=1/2}=\frac{\partial f_2(I_i,I_j,I_k)}{\partial I_k}\rvert_{I_i=1/2}&=&1\\
    \frac{\partial f_2(I_i,I_j,I_k)}{\partial I_i}\rvert_{I_i=1/2}&=&-2\, .
\end{eqnarray}
}
The latter expressions enter into the computation of \cref{eq:h_i} and thus we eventually obtain the following final form for the control terms $\tilde{h}^{(M)}$ and $h^{(M)}$ given in \cref{eq:h_iPartial2D,eq:h_iPartial}:
    \begin{eqnarray}\displaystyle
	\label{eq:FullControlTerm2D}
	\tilde{h}^{(M)}&:=&\frac{1}{2}\left(\frac{K_1}{N}\sum_{k=1}^MA_{ki}\cos(\theta_k-\theta_i)  \right) \times \left(-\frac{K_1}{N}\sum_{k=1}^MA_{ki}\frac{\cos(\theta_k-\theta_i)}{\omega_k-\omega_i} \right)\\
	&-&\frac{1}{2}
	\left(-\frac{K_1}{N}\sum_{k=1}^MA_{ki}\sin(\theta_k-\theta_i) \right) \times \left(-\frac{K_1}{N}\sum_{k=1}^MA_{ki}\frac{\sin(\theta_k-\theta_i)}{\omega_k-\omega_i} \right) \nonumber \\
	&+&\frac{1}{2}\sum_{j=1}^M\left\{\left(-
	\frac{K_1}{N}A_{ij}\cos(\theta_j-\theta_i) \right)\times\left(-\frac{K_1}{N}\sum_{k=1}^M A_{kj}\frac{\cos(\theta_k-\theta_j)}{\omega_k-\omega_j}  \right) \right. \nonumber \\
	&&\left.-\left(-\frac{K_1}{N}\sum_{k=1}^MA_{jk}\sin(\theta_k-\theta_j)\right)\times\left(\frac{K_1}{N}A_{ij}\frac{\sin(\theta_j-\theta_i)}{\omega_j-\omega_i} \right) \right\} \nonumber 
\end{eqnarray}
and

\begin{eqnarray}\displaystyle
	\label{eq:FullControlTerm}
	h^{(M)}&:=&\frac{1}{2}\left(\frac{K_1}{N}\sum_{k=1}^MA_{ki}\cos(\theta_k-\theta_i) + \frac{K_2}{N^2}\sum_{k,l=1}^MB_{ikl}\left[\cos(\theta_i+\theta_k-2\theta_l)+2\cos(\theta_k+\theta_l-2\theta_i)\right] \right) \\
	&&\times \left(-\frac{K_1}{N}\sum_{k=1}^MA_{ki}\frac{\cos(\theta_k-\theta_i)}{\omega_k-\omega_i} + \frac{K_2}{N^2}\sum_{k,l=1}^MB_{ikl}\left[\frac{\cos(\theta_i+\theta_k-2\theta_l)}{\omega_i+\omega_k-2\omega_l}-\frac{\cos(\theta_k+\theta_l-2\theta_i)}{\omega_k+\omega_l-2\omega_i}\right] \right) \nonumber \\
	&-&\frac{1}{2}
	\left(-\frac{K_1}{N}\sum_{k=1}^MA_{ki}\sin(\theta_k-\theta_i)+\frac{K_2}{N^2}\sum_{k,l=1}^MB_{ikl}\left[\sin(\theta_i+\theta_k-2\theta_l)-\sin(\theta_k+\theta_l-2\theta_i)\right] \right) \nonumber \\
	&&\times \left(-\frac{K_1}{N}\sum_{k=1}^MA_{ki}\frac{\sin(\theta_k-\theta_i)}{\omega_k-\omega_i} + \frac{K_2}{N^2}\sum_{k,l=1}^MB_{ikl}\left[-\frac{\sin(\theta_i+\theta_k-2\theta_l)}{\omega_i+\omega_k-2\omega_l}-2\frac{\sin(\theta_k+\theta_l-2\theta_i)}{\omega_k+\omega_l-2\omega_i}\right] \right) \nonumber \\
	&+&\frac{1}{2}\sum_{j=1}^M\left\{\left(\frac{K_2}{N^2}\sum_{k=1}^MB_{ijk}\left[\cos(\theta_j+\theta_i-2\theta_k)-2\cos(\theta_k+\theta_i-2\theta_j)-2\cos(\theta_j+\theta_k-2\theta_i)\right] -
	\frac{K_1}{N}A_{ij}\cos(\theta_j-\theta_i) \right)\right. \nonumber \\
	&&\times\left(-\frac{K_1}{N}\sum_{k=1}^M A_{kj}\frac{\cos(\theta_k-\theta_j)}{\omega_k-\omega_j} + \frac{K_2}{N^2}\sum_{k,l=1}^MB_{jkl}\left[\frac{\cos(\theta_j+\theta_k-2\theta_l)}{\omega_j+\omega_k-2\omega_l}-\frac{\cos(\theta_k+\theta_l-2\theta_j)}{\omega_k+\omega_l-2\omega_j}\right] \right) \nonumber \\
	&&-\left(-\frac{K_1}{N}\sum_{k=1}^MA_{jk}\sin(\theta_k-\theta_j)+\frac{K_2}{N^2}\sum_{k,l=1}^MB_{jkl}\left[\sin(\theta_j+\theta_k-2\theta_l)-\sin(\theta_k+\theta_l-2\theta_j)\right]\right) \nonumber \\
	&&\left.\times\left(\frac{K_1}{N}A_{ij}\frac{\sin(\theta_j-\theta_i)}{\omega_j-\omega_i} + \frac{K_2}{N^2}\sum_{k=1}^M\left[-\frac{\sin(\theta_j+\theta_i-2\theta_k)}{\omega_j+\omega_i-2\omega_k}+2\frac{\sin(\theta_k+\theta_i-2\theta_j)}{\omega_k+\omega_i-2\omega_j}+2\frac{\sin(\theta_j+\theta_k-2\theta_i)}{\omega_j+\omega_k-2\omega_i}\right]\right) \right\}\, . \nonumber 
\end{eqnarray}   

{\rev

By assuming the underlying hypergraph to have an all-to-all topology, the authors of \cite{asllani2018minimally} showed that the pairwise control term could be reduce to its dominant term given by~\cref{eq:pairwise_reduction}. By following the same ideas, and defining additional rescaled order parameters:
\begin{equation}
R_2e^{\iota\Psi_2}:=\frac{1}{N}\sum_je^{\iota2\theta_j}\, ,
\quad
\tilde{\tilde{R}}_ie^{\iota\tilde{\tilde{\Psi}}_i}:=\frac{1}{N^2}\sum_{j,k}\frac{e^{-\iota(\theta_j+\theta_k)}}{\omega_j+\omega_k-2\omega_i}
\quad
\text{and}
\quad \bar{R}_ie^{\iota\bar{\Psi}_i}:=\frac{1}{N^2}\sum_{j,k}\frac{e^{-\iota(2\theta_j-\theta_k)}}{2\omega_j-\omega_k-\omega_i}
\end{equation}
we get the following dominant term
\begin{eqnarray}
   h^{(N)}_i &\approx&-\frac{1}{2}K_1^2R\tilde{R}_i\cos(\Psi-\tilde{\Psi}_i)-\frac12 K_1K_2RR_2\tilde{R}_i\cos(\Psi_2-\Psi-\tilde{\Psi}_i)\\
&-&\frac12 2K_1K_2R^2\tilde{R}_i\cos(2\Psi-\tilde{\Psi}_i-\theta_i)
	-\frac12 K_1K_2R^2\tilde{R}_i\cos(2\Psi-2\theta_i)\cos(\tilde{\Psi}_i-\theta_i) \nonumber\\
	&-&\frac12 K_1K_2R\tilde{\tilde{R}}_i\cos(\Psi + \theta_i - \tilde{\tilde{\Psi}}_i)\nonumber\\
	&-&\frac12 K_1K_2RR_2\tilde{\tilde{R}}_i\cos(\Psi_2 + \theta_i - \Psi - \tilde{\tilde{\Psi}}_i)\nonumber\\
	&-&\frac12 K_2^2R^2\tilde{\tilde{R}}_i\cos(2\Psi -  \tilde{\tilde{\Psi}}_i) - \frac12 K_2R^2\cos(2\Psi-2\theta_i)K_2\tilde{\tilde{R}}_i\cos(\tilde{\tilde{\Psi}}_i-2\theta_i)\nonumber\\
	&-&\frac12 K_1K_2R\bar{R}_i\cos(\Psi - \bar{\Psi}_i) - \frac12 K_1K_2R\bar{R}_i\sin(\Psi-\theta_i)\sin(\bar{\Psi}_i-\theta_i)\nonumber\\
	&-&\frac12 K_1K_2RR_2\bar{R}_i\cos(\Psi_2 - \Psi - \bar{\Psi}_i) - \frac12 K_1K_2RR_2\bar{R}_i\sin(\Psi_2 - \Psi-\theta_i)\sin(\bar{\Psi}_i-\theta_i)\nonumber\\
	&-&K_2^2R^2\bar{R}_i\cos(2\Psi-\bar{\Psi}_i-\theta_i)\, . \nonumber
\end{eqnarray}

}

\section{Generalization to any interaction order}
\label{app:OrderD}

The proposed Hamiltonian function~\eqref{def:HOHam3} can be extended to allow for interactions of arbitrary order $d$. To do so one can add a term 
\begin{equation}
    V^{(d)}= \frac{K_d}{N^d}\sum_{i_0,\dots,i_d=1}^NA^{(d)}_{i_0,\dots,i_d}\sqrt[d+1]{I_{i_0}\dots I_{i_d}}\left(\sum_{j=0}^{d}\alpha_jI_{i_j}\right)\sin{\left(\sum_{j=0}^{d}\alpha_j\theta_{i_j}\right)}
\end{equation}
where $K_d$ is the $d$-order interaction strength, $\mathbf{A}^{(d)}$ is the $d$-{th} order adjacency tensor of the underlying hypergraph and $\{\alpha_0,\dots,\alpha_d\}$ are integer coefficients that sum to zero. 

The choice of the coefficients $\alpha$ defines the type of $d$-order interactions that will be represented in the higher-order system. A standard choice would be: $\alpha_0=-d$ and $\alpha_j=1$ for all $j>0$, but there can exist many others possibilities. As it happens with the term $V^{(2)}$ described above, several $d$-order interaction terms appear in the embedded HOKM dynamics. They correspond to all the different permutations of $\mathbf{\alpha}$ such that $\alpha_0<0$, and its opposite $-\mathbf{\alpha}$.

For example for $d=3$ we can consider 
\begin{equation}
    V^{(3)}=\frac{K_3}{N^3}\sum_{i,j,k,l=1}^NA^{(3)}_{i,j,k,l}\sqrt[4]{I_{i}I_jI_kI_l}\left(I_j+I_k+I_l-3I_i\right)\sin\left(\theta_j+\theta_k+\theta_l-3\theta_{i}\right),
\end{equation}
i.e., $\alpha_0=-3$ and $\alpha_1=\alpha_2=\alpha_3=1$. Then the HOKM resulting from the Hamiltonian function $H:=H_0+V^{(3)}$ restricted to the invariant torus $\mathbf{I}=\frac{1}{2}$ is 
\begin{equation}\label{eq:4Order1}
    \dot{\theta_i}=\omega_i+\frac{3}{2}\frac{K_3}{N^3}\sum_{j,k,l=1}^NA^{(3)}_{i,j,k,l}\left[\sin\left(\theta_j+\theta_k+\theta_l-3\theta_{i}\right)+\sin\left(3\theta_j-\theta_k-\theta_l-\theta_{i}\right)\right].
\end{equation}
If on the other hand one defines 
\begin{equation}
    V^{(3)}=\frac{K_3}{N^3}\sum_{i,j,k,l=1}^NA^{(3)}_{i,j,k,l}\sqrt[4]{I_{i}I_jI_kI_l}\left(I_k+I_l-I_j-I_i\right)\sin\left(\theta_k+\theta_l-\theta_j-\theta_{i}\right),
\end{equation}
i.e. $\alpha_0=\alpha_1=-1$ and $\alpha_2=\alpha_3=1$ then the resulting system is
\begin{equation}\label{eq:4Order2}
    \dot{\theta_i}=\omega_i+2\frac{K_3}{N^3}\sum_{j,k,l=1}^NA^{(3)}_{i,j,k,l}\sin(\theta_k+\theta_l-\theta_j-\theta_i).
\end{equation}

{\rev Let us finally remark that }the multiplying coefficients $\frac{3}{2}$, resp. $2$, appearing in \cref{eq:4Order1} and \cref{eq:4Order2} will be replaced by $1$ if one sets $\mathbf{I}=\frac{1}{3}$, resp. $\mathbf{I}=\frac{1}{4}$, instead of $\mathbf{I}=\frac{1}{2}$. The theory developed stays valid independently of the constant choice. If one has to consider order 2 and 3 together, for example, a re-scaling of $K_3$ can be made to make the multiplicative constants disappear.

\section{Larger all-to-all hypergraph}
\label{app:N=100}

In this section is we present some results for larger hypergraph where $N=100$ nodes are connected in an all-to-all fashion. Here again we set $(K_1,K_2)\in[0,2]$, $\theta_0\sim U([0,0.3])$ and $\mathbf{\omega}\sim U([0,1])$. 

\begin{figure}[hbt]
    \centering
    \includegraphics[width=0.8\linewidth]{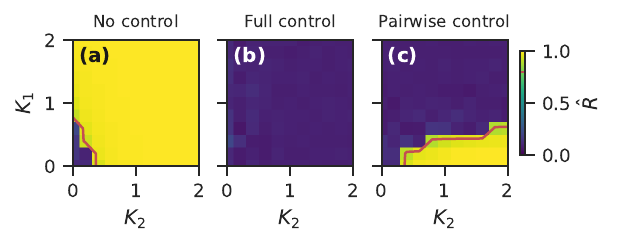}
    \caption{Desynchronizing by controlling higher-order interactions in an all-to-all hypergraph.  This figure is similar to \cref{Fig:All-to-all50}c-e but with $N=100$.  
    The average order parameter is shown over parameter space in three cases: (a) no control, (b) full control, and (c) pairwise control. 
    The boundary of the synchronized region is highlighted by a level curve ($\hat{R}=0.8$) in red. }
    \label{Fig:All-to-all100}
\end{figure}

In Fig.~\ref{Fig:All-to-all100} we compare the values of $\hat{R}$ reached by the uncontrolled and controlled systems by using terms $h^{(N)}$ or $\tilde{h}^{(N)}$. The results are analogous to the case $N=50$: $h^{(N)}$ enables to desyncrhonize the system for all tested $(K_1,K_2)$ values ($\hat{R}$ reaches the maximal value of $0.21$) and $\tilde{h}^{(N)}$ can desynchronize the system in all cases but when $K_1$ is small compared to $K_2$ (see the region bounded by the red curve in Fig. \ref{Fig:All-to-all100}(c)).

\begin{figure}[t!]
    \centering
    \includegraphics[width=0.7\linewidth]{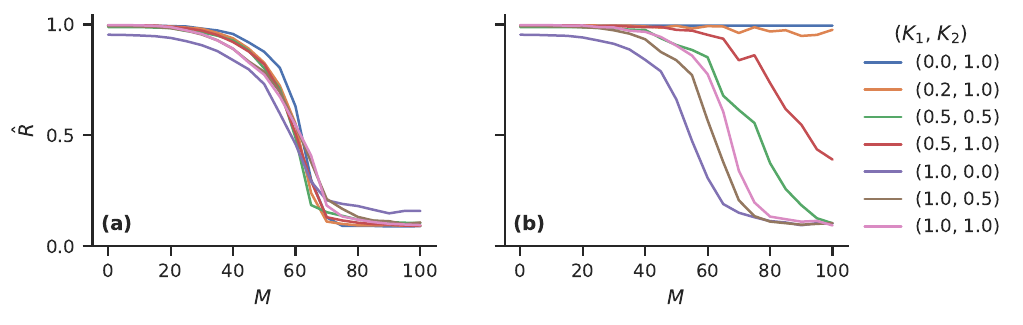}
    \caption{Average Kuramoto order parameter, $\hat{R}$, as a function of the number $M$ of controlled nodes in a complete hypergraph with $N=100$ nodes. In the panel (a) we show the results obtained by using $h^{(M)}$ as the control term; in (b) we report $\hat{R}$ resulting from the use of the control limited to the second order interactions only, i.e. $\tilde{h}^{(M)}$. Each curve has been obtained by fixing the coupling strengths $(K_1,K_2)$, the used values are reported in the legends, and it is the average of $100$ independent numerical simulations corresponding to different samples for $\pmb{\omega}\sim U([0,1])$.}
    \label{Fig:All-to-all100:Partial}
\end{figure}

Fig.~\ref{Fig:All-to-all100:Partial} shows $\hat{R}$ in function of $M$ for $h^{(M)}$ and $\tilde{h}^{(M)}$. Here again one can see in both cases a sharp decrease of $\hat{R}$ for $M\in[2N/5,3N/5]$ after what it roughly reaches its minimal value. The only exceptions concern once again $\tilde{h}^{(M)}$ when $K_1$ is ``small'' and $K_2$ is ``large'' (yellow region in the bottom right of Fig.~\ref{Fig:All-to-all100}(c)). 

\section{Basins of attractions and other states}
\label{App:Basins}

Here, we explore the possible equilibria of the all-to-all case, and the relative size of their basins of attraction. First, we observed only one state other than full synchronization and incoherence: 2-cluster states, where oscillators are divided into two clusters separated by a distance of $\pi$ (see e.g.~\cite{skardal2019abrupt, zhang2023deeper}). To compute the relative basin size of each state, we simulated 100 random initial conditions and automatically identified the equilibrium they reached. The relative basin size of a state is then simply computed as the number of random initial conditions that reach it, divided by the total number of initial conditions. Note that a 2-cluster state can be more or less balanced depending on the relative sizes of its clusters, that is, the number of oscillators in each cluster. 

\begin{figure}[t!]
    \centering
    \includegraphics[width=0.75\linewidth]{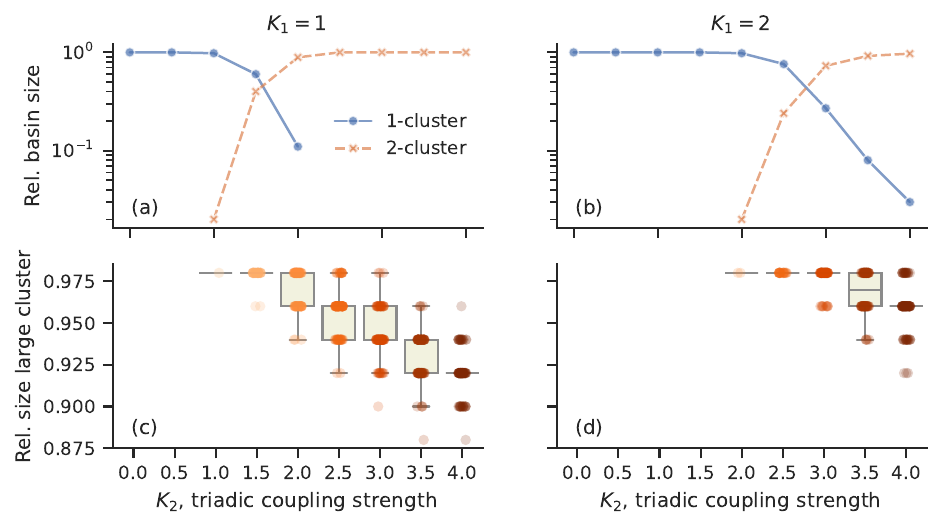}
    \caption{Relative basin sizes for the all-to-all coupling scheme, for (a, c) $K_1=1$ and (b, d) $K_1=2$, and a range of $K_2$ values. For both cases, we show (a,b) the relative basin sizes, and (c,d) the relative size of the larger cluster in 2-cluster states. The 2-cluster states are very unbalanced: one cluster contains about 90-100\% of the nodes. Hence, they are close to full synchronization.}
    \label{fig:basin_sizes}
\end{figure}

\cref{fig:basin_sizes} shows these basin sizes for $K_1=1$ and  $K_1=2$. First, we see that, for weak enough triadic coupling $K_2$, full synchronization (1-cluster) is the only attractor. Then as $K_2$ increases, 2-cluster states appear and quickly take over the phase space, with the basin of full synchronization shrinking dramatically. Note that for the stronger pairwise coupling $K_1$, that switch occurs at a lager $K_2$. Finally, we see that the 2-cluster states are interestingly very unbalanced: one cluster contains above 90\% of the nodes. This indicates that, even when $K_2$ is large enough and many initial conditions reach a 2-cluster state, that state is still very close to being full synchronization and this should not affect the control method.

\section{Analysis of the control cost}
\label{app:Cost}

In the previous sections we have observed that the various control strategies can have different outcomes, in this section we consider the cost we can associate to each one of them; more precisely we study the energy required for the control term to keep the system far from the synchronization manifold. For a given time interval $[0,T]$, this translates into the following formula
\begin{eqnarray}\label{eq:defCost}
    C^{(M)}:=\frac{1}{TM}\sum_{i=1}^M\int_0^{T}|h_i^{(M)}(t)|dt&; &
    \tilde{C}^{(M)}:=\frac{1}{TM}\sum_{i=1}^M\int_0^{T}|\tilde{h}_i^{(M)}(t)|dt.
\end{eqnarray}

\begin{figure}[ht]
    \centering
    \includegraphics[width=0.8\linewidth]{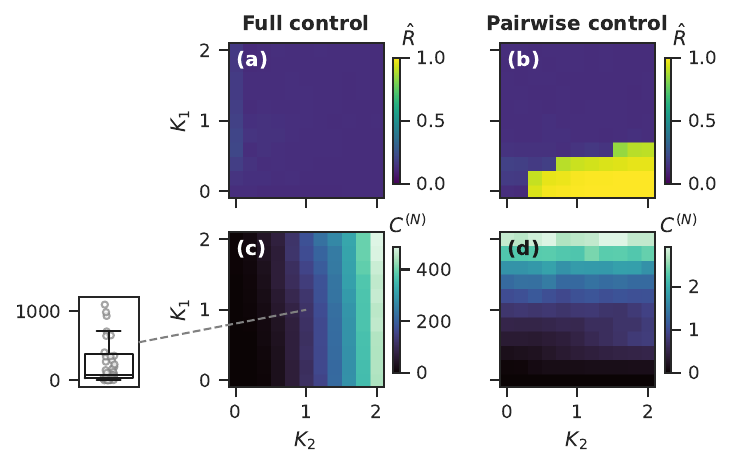}
    \caption{\rev \textbf{Comparison of cost.} (a,b) Average order parameter and associated cost (c,d) in parameter space, for (a,c) the full control and (b,d) the pairwise control. 
    The structure is an all-to-all hypergraph with $N=50$. The color maps display the median values obtained from 50 simulations with sets of initial conditions and natural frequencies following $\theta_0\sim U[0,0.3]$ and $\omega\sim U[0,1]$. The integration \eqref{eq:defCost} was then computed with trapeze method. 
    The inset in (c) show the distribution of values for $K_1=K_2=1$.}
    \label{Fig:ControlCost}
\end{figure}

In Fig.~\ref{Fig:ControlCost} we show the cost functions $\tilde{C}^{(N)}$ and $C^{(N)}$ as a function of the coupling strengths $K_1$ and $K_2$, for the choice $T=40$. It is clear that the control $h^{(N)}$ is always more costly than  $\tilde{h}^{(N)}$ by several orders of magnitude. This can be explained as the terms involved in the former contain double and triple sums and therefore increasing cost. This result suggests that the use of $\tilde{h}^{(N)}$ is preferable once $K_1$ is not too small compared to $K_2$.

In addition we can observe that $C^{(N)}$ increases monotonically with $K_2$ but does not vary significantly in function of $K_1$. This fact strengthen the last observation about the dominance of the higher-order terms in the control cost. On the other hand, $\tilde{C}^{(N)}$ increases with $K_1$ but does not vary with $K_2$. This indicates that the control effort that is necessary to desynchronize the system by using $\tilde{h}^{(N)}$ does not increase even if $K_2$ increases and, as a consequence, the attraction of the synchronization state on the close trajectories increases.

Finally, we can note that the different cost values can reach quite diverse values, especially when $C^{(N)}$. In the inset of Fig. \ref{Fig:ControlCost}c we present this dispersion by using a boxplot, one can observe that the median value is much smaller that the maximum value and quite close to the minimum one. There are indeed some configurations of $\pmb{\omega}$ and $\pmb{\theta}_0$ giving rise to higher values of control cost. In particular, we found one of the tested configurations (that have been removed from the averages computations in Fig. \ref{Fig:ControlCost} because of its outlier behaviour) that gave rise to $C^{(N)}$ values of order $10^7$. There are multiple possible explanations to this phenomenon. There could be some region of the basin of attraction of the synchronized state, from which one needs a lot more energy to be extracted out than other regions. More likely, as the natural frequencies are directly present in the control definition, the particular values certainly impact significantly the control cost. Those points would deserve deeper investigations in further work.

\end{document}